\documentclass[12pt,a4paper,twoside]{amsart}

\usepackage{geometry}
\usepackage{fancyhdr}
\usepackage{txfonts}
\usepackage{bbm}
\usepackage{threeparttable}
\usepackage{mathrsfs}
\usepackage{amsfonts}
\usepackage{amsmath}
\usepackage{amssymb}
\usepackage{wasysym}
\usepackage{psfrag}
\usepackage{graphics,epsfig,amsmath}  
\usepackage{comment}

\usepackage{xcolor}
\usepackage{enumerate}
\usepackage{hyperref}
\usepackage[bbgreekl]{mathbbol}
\usepackage{bm}

\DeclareSymbolFontAlphabet{\mathbb}{AMSb}
\DeclareSymbolFontAlphabet{\mathbbl}{bbold}

\newtheorem{theorem}{Theorem}[section]

\newtheorem{proposition}{Proposition}

\newtheorem{lemma}{Lemma}[section]

\newtheorem{remark}{Remark}[section]

\newtheorem{definition}{Definition}

\newcommand{\ZZ}{\mathbb{Z}}

\newcommand{\RR}{\mathbb{R}}

\newcommand{\TT}{\mathbb{T}}

\newcommand{\scr}{\mathscr}

\allowdisplaybreaks

\numberwithin{equation}{section}

\title[Quasi-periodic equilibrium]{KAM theory for quasi-periodic equilibrium of generalized Frenkel-Kontorova model in $\mathbb{Z}^{d}$}

\author[X. He]{Xiaolong He}
\address{
School of mathematics\\
Hangzhou Normal University\\
Hangzhou, 311121, P.R. China}

\email{xlhe@hznu.edu.cn}

\author[Z. You]{Zuhong You}
\address{
School of Mathematical Science\\
Fudan University\\
220 Handan Road, Shanghai, 200433, People's Republic of China}

\email{23110180049@m.fudan.edu.cn}

\begin{document}

\vspace{0.1in}

\begin{abstract}
We extend the KAM theory for  one dimensional  Frenkel-Kontorova
model to  multidimensional case in quasi-periodic media.
By using the quasi-Newton method, we reduce the linearized operator with quasi-periodic
coefficient in
the Nash-Moser iteration
to the composition of two diagonal operators and a near-identity operator,
such that the small divisor problem can be overcome by imposing the simultaneous Diophantine condition.
In contrast,
we also discuss the application of Craig-Wayne-Bourgain method
to the quasi-periodic equilibrium problem.

\end{abstract}


\maketitle

\section{Introduction}

In this paper, we study the quasi-periodic equilibrium of the generalized 
Frenkel-Kontorova (F-K) model 
\begin{equation}\label{FK}
(\mathbf{\Delta}_{d} u)_{n}+ \epsilon V'(u_{n})=0, \quad n\in\mathbb{Z}^{d},
\end{equation}
where $V$ is a quasi-periodic function with frequency $\alpha\in\mathbb{R}^{m}$ and the discrete Laplacian $\mathbf{\Delta}_{d}$ is given by
\begin{equation*}
(\mathbf{\Delta}_{d} u)_{n}= \frac{1}{2d} \sum_{j=1}^{d} u_{n+e_{j}}+u_{n-e_{j}}-2u_{n},
\end{equation*}
with the standard basis  $\{e_{1},\cdots, e_{d}\}$ in $\RR^{d}$.

When $d=1$, the KAM theory for the  quasi-periodic equilibrium of  one dimensional F-K model is established in 
\cite{SdlL12b, SdlL12a}. 
The goal of this paper is to generalize the results in \cite{SdlL12b, SdlL12a}
to the case of $d\geq 2$ as follows.

\begin{theorem}\label{KAM}
Let $\alpha\in\mathbb{R}^{m}$, $\omega\in\RR^{d}$ and $\tau>2m+1$.
Assume that 
\begin{enumerate}
\item [(i)] the function
$V$ is a real analytic and  quasi-periodic  with frequency $\alpha\in\mathbb{R}^{m}$;
\item [(ii)] the collection of vectors $\omega_{1}\alpha, \cdots, \omega_{d}\alpha$ satisfies
the simultaneous Diophantine condition 
\begin{equation*}
\max_{1\leq j\leq d} \inf_{n\in \mathbb{Z}} |\omega_{j}\langle k, \alpha\rangle-n|>\frac{\gamma}{|k|^{\tau}},\quad 
\forall~k\in\mathbb{Z}^{m}\setminus\{0\},
\end{equation*}
where $|k|= \sum_{j=1}^{m} |k_{j}|$ and $\langle k,\alpha\rangle=\sum_{j=1}^{m} k_{j} \alpha_{j}$.
\end{enumerate}
Then there exists $\epsilon_{\ast}=\epsilon_{\ast}(d, m, \alpha, \gamma, \tau, V)>0$ such that if $0<\epsilon< \epsilon_{\ast}$, there exists a real analytic function $h: \TT^{m}\rightarrow\RR$ with $\int_{\TT^{m}} h(\theta)~d\theta=0$ such that
\begin{equation}\label{config 1}
u_{n}= \langle n, \omega\rangle+ h( \langle n, \omega\rangle \alpha)
\end{equation}
is a solution of \eqref{FK}.
\end{theorem}

The proof of Theorem \ref{KAM} is delayed to Section \ref{proof KAM}. Due to the small divisor problem,
the function $h$ loses some analyticity domain in contrast with  that of  the hull function of $V$.
Furthermore, in Section \ref{sect SDC}, it is shown that the set of $\omega$ such that
the simultaneous Diophantine condition occupies full Lebesgue measure in
$\RR^{d}$ provided $\alpha$ is Diophantine and the exponent $\tau$ in \eqref{SDC} is large.

\begin{remark}
Note that the quasi-periodic states of \eqref{FK} are not unique. If we can find a function $h: \TT^{m}\rightarrow \RR$
such that $u_{n}$ given by \eqref{config 1} is a solution, then one readily sees that for any $s\in\RR$,
\begin{equation*}
u_{n}= \langle n, \omega\rangle+s+h(\langle n, \omega\rangle\alpha+s\alpha)
\end{equation*} 
is also a solution of \eqref{FK}.
For this reason, we look for the solution $h$ of \eqref{FK} such that its average vanishes
in Theorem \ref{KAM}.

\end{remark}

The equilibrium configuration of the F-K models has attracted lots of attention in recent years.
We refer to \cite{dlL08,dlLV07} for the physical description of the models and the dynamical interpretation of the
equilibrium state solutions.
There are usually two approaches for constructing solutions of the equilibrium equation: one is
the KAM theory and the other is the variational theory.
The KAM theory for periodic and quasi-periodic equilibrium for F-K model could be found in
the papers \cite{dlL08, dlLSZ16, dlLSZ17, SdlL12b, SdlL12a} and the references therein.
Recently, the KAM theory for almost-periodic equilibrium for one dimensional F-K model is
also established in \cite{AdlLSWY26}.

In contrast,
 the variational method usually does not need Diophantine
conditions on the frequencies and has moderate regularity assumptions on the interaction, while
it requires convexity assumptions that are not needed in KAM theory.
We refer to \cite{CdlL98, dlLV07, GGP06, SdlL17,Tre19} and the references therein for more details.

In this paper, we  study the analytic quasi-periodic equilibrium of \eqref{FK} by the KAM theory. 
In each step of the Nash-Moser iteration, we apply the quasi-Newton method
to reduce the linearized operator with quasi-periodic coefficient into
the composition of two diagonal operators and a near identity operator,
such that the small divisor problem appears only in the diagonal part.
As a result, we only need to solve the homological equation with constant
coefficient by using the simultaneous Diophantine condition.

It is  possible 
to generalize our results to the case of non-nearest neighborhood or to the case of Gevrey or finitely differentiable  regularity.  Moreover, one can 
also formulate the KAM theorem in an a-posteriori form as \cite{dlL08, SdlL12a}. We do not 
pursue them here.  
Instead, we would like to compare the classical KAM method developed here with
the Craig-Wayne-Bourgain method (see Section \ref{sect dis}),  which is a well-known approach for its advantage of solving linear equations with
variable coefficients.

\subsection{Formulation of the functional equation}\label{sect FE}
To save notations, we still denote  by $V$ the associated hull function and write $V(x)= V(\alpha x)$ with   $V: \mathbb{T}^{m}=\mathbb{R}^{m}/\mathbb{Z}^{m}\rightarrow \mathbb{R}$. It turns out that
$V'(x)= \partial_{\alpha} V(\alpha x)=\sum_{j=1}^{m} \frac{\partial V}{\partial \theta_{j}}(\alpha x) \cdot \alpha_{j}.
$
Then the original model
\eqref{FK} can be written as
\begin{equation}\label{FK'}
(\mathbf{\Delta}_{d} u)_{n}+ \epsilon \partial_{\alpha}V( \alpha u_{n})=0, \quad n\in\mathbb{Z}^{d}.
\end{equation}

We are looking for a function $h: \mathbb{T}^{m}\rightarrow \mathbb{R}$ with the average
$\langle h\rangle$ of $h$ being zero such that $u_{n}$ given by \eqref{config 1}
is a solution of \eqref{FK'}.
Substituting \eqref{config 1} into \eqref{FK'} yields
\begin{equation*}
\begin{aligned}
\frac{1}{2d} \sum_{j=1}^{d} \left[h(\langle n,\omega\rangle \alpha+\omega_{j}\alpha )+h(\langle n,\omega\rangle \alpha-\omega_{j}\alpha )-2h(\langle n,\omega\rangle \alpha)\right]
+\epsilon \partial_{\alpha} V(\langle n,\omega\rangle \alpha+ \alpha h(\langle n,\omega\rangle \alpha) )=0.
\end{aligned}
\end{equation*}
Assuming $\alpha$ and $\omega$ are irrational vectors, it suffices to solve $h$ from
\begin{equation}\label{FE}
\frac{1}{2d}\sum_{j=1}^{d}\left[ h(\theta+\omega_{j}\alpha )+h(\theta-\omega_{j}\alpha )-2h(\theta )\right]
+\epsilon \partial_{\alpha} V(\theta+ \alpha h(\theta) )=0,\quad \theta\in \mathbb{T}^{m}.
\end{equation}

For simplicity, we define  the translation operator $(T_{\beta} h)(\theta)=h(\theta+\beta) $ for  $\beta\in\mathbb{T}^{m}$ and the difference operator $\Delta$ as follows
\begin{equation}
\label{def Delta}
\Delta h=  \frac{1}{2d} \sum_{j=1}^{d}\left[ T_{\omega_{j} \alpha} h+T_{-\omega_{j}\alpha}h -2 h\right].
\end{equation}
Following \cite{SdlL12a}, we define the functional $\scr{F}$ with extra one dimensional parameter
\begin{equation}\label{def F}
\scr{F}[h, \lambda](\theta)=\Delta h (\theta)+\epsilon U(\theta+\alpha h(\theta))+\lambda,
\end{equation}
and study the functional equation
\begin{equation}\label{FE'}
\scr{F}[h, \lambda]=0.
\end{equation}

It is shown in \cite[Lemma 23]{SdlL12a} that whenever \eqref{FE'} has a variational structure (i.e., $U=\partial_{\alpha} V$), the solution of \eqref{FE'} leads
 to $\lambda=0$. See Lemma \ref{vanishing lm 1}  for our  case of $d\geq 2$. 
The extra parameter $\lambda$ allows a more efficient iterative procedure
as introduced in \cite{Mos67}.
See \cite[Section 1.2]{SdlL12a} for some explanations on the external forces.

We shall apply the Nash-Moser iteration method to solve \eqref{FE'}. Let us briefly explain 
the main idea and the advantage of the method developed in \cite{dlL08}, which can also
be adapted to the high dimensional problem in this paper.

Given  an approximate solution $(h, \lambda)$ of \eqref{FE'}, the Newton equation for the correction 
$(\delta_{h}, \delta_{\lambda})$ is given by
\begin{equation}\label{Newton 1}
\Delta \delta_{h} (\theta)+\epsilon \partial_{\alpha} U(\theta+ \alpha h(\theta))\cdot \delta_{h}(\theta)+\delta_{\lambda}=- \scr{F}[h, \lambda](\theta),
\end{equation}
which contains the variable coefficient $\partial_{\alpha} U(\theta+ \alpha h(\theta))$.
One possible way to solve \eqref{Newton 1} is to apply the well-developed Green's function estimate method,
which seems technically involved (see Section \ref{sect dis} for further discussion).

A key observation  is a geometry identity (see \cite[Section 5.3]{dlL08}) that enables us to reduce
\eqref{Newton 1} into a constant coefficient equation up to a quadratic error.
This is exactly the case of $d=1$ developed in \cite{SdlL12a}. However, for $d\geq 2$, 
the transformed linearized operator, denoted by $L_{h}$, still contains the variable coefficient
but can be written as 
\begin{equation}\label{Lh trans}
L_{h}= S_{w} (id+K_{h}) S_{w},
\end{equation}
in which $S_{w}$ is a diagonal operator.

The novelty is that  small divisor problem appears only in constructing the inverse of  $S_{w}$,
which can be implemented by the simultaneous Diophantine condition (see Lemma \ref{Sw -1}).
In contrast, the operator $id+K_{h}$ has a bounded inverse by simply using the Neumann series.
The idea of the factorization as \eqref{Lh trans} is very close
to that in \cite{Koz83} on the reducibility problem of elliptic differential operators\footnote{We thank Prof. Rafael de la Llave for bringing our attention to \cite{Koz83} and sharing the manuscript.}.

We emphasize that, in \cite{dlL08}, similar transformed operator like \eqref{Lh trans} also
appears in the periodic equilibrium problem for 1-D statistical mechanics models with
non-nearest  neighborhood. 

\subsection{Organization of the paper.}
The rest of the paper is organized as follows. In Section \ref{sect pre}, we collect the function spaces and some useful properties that are used frequently later.
In Section \ref{proof KAM}, we prove Theorem \ref{KAM} by applying the 
Nash-Moser iteration and the quasi-Newton method.
In Section \ref{sect dis}, we discuss briefly the idea of applying the Craig-Wayne-Bourgain method to the quasi-periodic equilibrium problem, and compare it with
the classical KAM method developed in the present paper.

\section{Preliminary}\label{sect pre}
In this section, we collect the function spaces and the non-resonant conditions for the frequencies. Then we solve the homological equation obtained from the quasi-Newton method. Finally, we prove the vanishing lemma to show that the additional parameter $\lambda$
in \eqref{FE'} can
be eliminated.
Following the conventional notations in KAM theory, we denote in the sequel
various constants by the letter $C$ with some subscripts indicating the dependence
on the given quantities, which might be different from line to line. These constants
could be made explicit from the context.

\subsection{Function space}

Let $\mathbb{T}^{m}=\RR^{m}/\ZZ^{m}$.
Given a periodic function $f:\TT^{m}\rightarrow \RR$, we expand it into Fourier series
\begin{equation*}
f(\theta)= \sum_{k\in\ZZ^{m}} \widehat{f}(k) e^{2\pi i \langle k,\theta\rangle},\quad i=\sqrt{-1},
\end{equation*}
where $\widehat{f}(k)$ is the Fourier coefficient.
The average of $f$ is denoted by
\begin{equation*}
\langle f\rangle= \int_{\TT^{m}} f(\theta)~d\theta= \widehat{f}(0).
\end{equation*}
Given $\alpha\in\RR^{m}$, we denote by $\partial_{\alpha} f$ the directional derivative
\begin{equation*}
\partial_{\alpha}f(\theta)= Df(\theta)\cdot \alpha= \sum_{j=1}^{m}\frac{\partial f(\theta)}{\partial \theta_{j}}\cdot \alpha_{j}.
\end{equation*}

Define the weighted $\ell^{1}$-norm
\begin{equation*}
\|f\|_{\rho}=\sum_{k\in\ZZ^{m}} |\widehat{f}(k)|  \cdot  e^{2\pi \rho |k|}.
\end{equation*}
Under the above norm, the set of periodic functions
 \begin{equation*}
 \scr{A}_{\rho}=\left\{f: \TT^{m}\rightarrow \RR~|~\|f\|_{\rho}<\infty\right\}
 \end{equation*}
 is a Banach space. 
  We also denote the closed subspace of $\scr{A}_{\rho}$ such that the functions  have  vanishing average  by
 \begin{equation*}
 \scr{A}_{\rho}^{0}=\left\{f: \TT^{m}\rightarrow \RR~|~\langle f\rangle=0,\quad \|f\|_{\rho}<\infty\right\}.
 \end{equation*}
 Let $\scr{L}(\scr{A}_{\rho}^{0},\scr{A}_{\rho'}^{0})$ be the Banach space of the bounded
 linear operators from  $\scr{A}_{\rho}^{0}$ to $\scr{A}_{\rho'}^{0}$, whose induced operator norm is
 denoted by $\|\cdot \|_{\scr{L}(\scr{A}_{\rho}^{0},\scr{A}_{\rho'}^{0})}$. 

 We collect some properties of functions in $\scr{A}_{\rho}$.
 
 \begin{lemma}
 Let $0<\sigma<\rho$ and $\alpha\in \mathbb{R}^{m}$.
 \begin{enumerate}[(i)]
 \item If $f, g\in \scr{A}_{\rho}$, then $\|f\cdot g\|_{\rho}\leq \|f\|_{\rho}\cdot \|g\|_{\rho}$.
 \item If $f\in \scr{A}_{\rho}$, then $\|\partial_{\alpha} f\|_{\rho-\sigma}\leq C_{\alpha, m}  \sigma^{-1}\|f\|_{\rho}$.
 \end{enumerate}
 \end{lemma}
 
 \noindent\textbf{Proof.~} (i) The Fourier coefficient $\widehat{fg}(k)$ is given by the convolution product
 \begin{equation*}
 \widehat{fg}(k)=\sum_{k'\in\ZZ^{m}} \widehat{f}(k-k')\cdot \widehat{g}(k').
 \end{equation*}
 Then we have
 \begin{equation*}
 \begin{aligned}
 \|fg\|_{\rho}=& \sum_{k\in\ZZ^{m}} \left|\sum_{k'\in\ZZ^{m}} \widehat{f}(k-k')\cdot \widehat{g}(k')\right|\cdot e^{2\pi \rho |k|}\\
 \leq & \sum_{k'\in\ZZ^{m}}   |\widehat{g}(k')| e^{2\pi \rho |k'|}\cdot \sum_{k\in\ZZ^{m}}  |\widehat{f}(k-k')|\cdot e^{2\pi \rho |k-k'|}
 \leq \|f\|_{\rho}\cdot \|g\|_{\rho}.
 \end{aligned}
 \end{equation*}
  
(ii) Observing that $\widehat{\partial_{\alpha}f}(k)= 2\pi i \langle k,\alpha\rangle \widehat{f}(k)$, we have
\begin{equation*}
\begin{aligned}
\|\partial_{\alpha}f\|_{\rho-\sigma}=& \sum_{k\in\ZZ^{m}} |\langle k,\alpha\rangle|\cdot | \widehat{f}(k)|\cdot e^{2\pi (\rho-\sigma)|k|}\\
\leq & |\alpha|\cdot \|f\|_{\rho}\cdot \sum_{k\in\ZZ^{m}} |k| e^{-2\pi \sigma |k|}\leq C_{\alpha, m} \sigma^{-1} \|f\|_{\rho},
\end{aligned}
\end{equation*}
where the constant $C_{\alpha, m}$ depends only on $\max_{1\leq j\leq m}|\alpha_{j}|$ and $m$.
\qed

\subsection{Simultaneous Diophantine condition}\label{sect SDC}
As the generalization of the well-known Diophantine condition, we define the simultaneous 
Diophantine condition as follows.
\begin{definition}
A collection of vectors $\bm{v}_{1}, \bm{v}_{2}, \cdots, \bm{v}_{d} \in \mathbb{R}^{m}$ is
 \textbf{simultaneous Diophantine}  of type $(\gamma, \tau)$ if
\begin{equation}\label{SDC}
\max_{1\leq j\leq d} \inf_{n\in \mathbb{Z}} |\langle k, \bm{v}_{j}\rangle-n|>\frac{\gamma}{|k|^{\tau}},\quad 
\forall~k\in\mathbb{Z}^{m}\setminus\{0\},
\end{equation}
where $|k|= \sum_{j=1}^{m} |k_{j}|$.
\end{definition}

Obviously, if one of the vectors satisfies the well-known Diophantine condition, then
the simultaneous Diophantine condition holds. In what follows, we shall
encounter with the vectors $\bm{v}_{j}=\omega_{j}\alpha$ ($1\leq j\leq d$)
with the same direction $\alpha$.

We recall some results in \cite{SdlL12a} for
completeness.
When $d=1$, it is shown that if $\alpha\in\RR^{m}$ satisfies
\begin{equation}\label{DC alpha}
|\langle k, \alpha\rangle|> \mu |k|^{-\nu},\quad \forall~k\in\ZZ^{m}\setminus\{0\},
\end{equation}
then the set of $\omega\in \mathbb{R}$ such that $\omega\alpha$ is Diophantine as
\begin{equation*}
\inf_{n\in \mathbb{Z}} |\langle k, \omega\alpha\rangle-n|>\frac{\gamma}{|k|^{\tau}},\quad 
\forall~k\in\mathbb{Z}^{m}\setminus\{0\},
\end{equation*}
for some $\gamma>0$ occupies full Lebesgue measure in $\mathbb{R}$ if $\tau> m+\nu$.

As a result of \cite[Lemma 10]{SdlL12a},
in our case of  $d\geq 2$, one readily sees that the set of $\omega\in \mathbb{R}^{d}$ such that the collection of vectors
$\omega_{1}\alpha, \cdots, \omega_{d}\alpha$ satisfies the simultaneous Diophantine
condition \eqref{SDC} occupies the full Lebesgue measure in $\mathbb{R}^{d}$ if $\tau>\nu+m$ and  \eqref{DC alpha} 
 holds.

\begin{remark}
When $m=1$ and $d\geq 2$, the simultaneous Diophantine condition \eqref{SDC} appears in \cite{Mos90} to establish
the simultaneous conjugation of commuting circle diffeomorphisms. For $m\geq 2$ and $d\geq 2$,
\eqref{SDC} can be used to establish the simultaneous conjugation of commuting
foliation preserving torus maps \cite{He23}, which are related to the analyticity problem
of quasi-periodic solutions for state-dependent delay differential equations.
\end{remark}

\subsection{Homological equation}
For any $k\in\ZZ^{m}$, we define the weight function
\begin{equation}\label{wei}
w(k)=\left\{\frac{1}{2d}\sum_{j=1}^{d} (e^{ 2\pi i \omega_{j}\langle k, \alpha\rangle}-1)(e^{- 2\pi i  \omega_{j}\langle k, \alpha\rangle}-1)\right\}^{1/2}=\left\{ \frac{2}{d} \sum_{j=1}^{d}
\sin^{2}\left(\pi \omega_{j}\langle k, \alpha\rangle\right)\right\}^{1/2}\geq 0.
\end{equation}
Given $f(\theta)= \sum_{k\in \ZZ^{m}}\widehat{f}(k) e^{2\pi i\langle k,\theta\rangle}$, we
define the linear diagonal operator $S_{w}$ by
\begin{equation*}
(S_{w} f)(\theta)= \sum_{k\in\ZZ^{m}} w(k) \widehat{f}(k) e^{2\pi i \langle k, \theta\rangle}.
\end{equation*}

The following lemma is denoted to establishing the inverse operator of $S_{w}$ as well as its estimate.
\begin{lemma}\label{Sw -1}
Assume that the  collection of vectors $\omega_{1}\alpha, \cdots, \omega_{d}\alpha$ satisfies the simultaneous Diophantine condition
\eqref{SDC}. Then  for any $0<\sigma<\rho$, the inverse operator $S_{w}^{-1}: \scr{A}^{0}_{\rho}\rightarrow \scr{A}^{0}_{\rho-\sigma}$
exists and satisfies
\begin{equation*}
\|S_{w}^{-1} g\|_{\rho-\sigma}\leq C_{d,m,\tau}~\gamma^{-1}~ \sigma^{-\tau-m}~ \|g\|_{\rho}.
\end{equation*}
\end{lemma}
\noindent\textbf{Proof.~} 
Consider the equation $S_{w} h=g$ with $g\in \scr{A}_{\rho}^{0}$. By matching the
Fourier coefficients, we obtain a unique solution
\begin{equation*}
h(\theta)= \sum_{0\neq k\in\ZZ^{m}} \frac{1}{w(k)} \widehat{g}(k) e^{2\pi i \langle k,\theta\rangle}
\end{equation*}
such that $\langle h\rangle=0$.
If $\omega_{1}\alpha,\cdots, \omega_{d}\alpha$ fulfills \eqref{SDC},
we have
\begin{equation*}
w(k)\geq \frac{1}{\sqrt{d}} \max_{1\leq j\leq d} \left|\sin (\pi \omega_{j} \langle k,\alpha\rangle)\right|
\geq \frac{2}{\sqrt{d}} \max_{1\leq j\leq d} \inf_{n\in\ZZ}|\omega_{j} \langle k,\alpha\rangle-n|\geq
\frac{2}{\sqrt{d}}\cdot \frac{\gamma}{ |k|^{\tau}}.
\end{equation*}
As a result, there is
\begin{equation*}
\|h\|_{\rho-\sigma}\leq \frac{\sqrt{d}}{2\gamma} \sum_{k\in\ZZ^{m}} |k|^{\tau}\cdot |\widehat{g}(k)|\cdot e^{2\pi(\rho-\sigma)|k|}\leq \frac{\sqrt{d}}{2\gamma}\|g\|_{\rho}\sum_{k\in\ZZ^{m}} |k|^{\tau} e^{-2\pi \sigma |k|}
\leq  C_{d,m,\tau} \gamma^{-1} \sigma^{-\tau-m} \|g\|_{\rho}.
\end{equation*}
This completes the proof. \qed

\begin{remark}
When $d=1$, the operator $S_{w}$ refers to the classical homological equation 
\begin{equation*}
h(\theta+\omega_{j}\alpha)-h(\theta)=g(\theta), \quad \omega_{j}\in\RR, \alpha\in\RR^{m},~\theta\in\TT^{m},
\end{equation*}
which can be solved under the assumption that $\omega_{j}\alpha$ is Diophantine.
See \cite{SdlL12a} for more details.
\end{remark}

\subsection{Vanishing lemma}\label{vanishing lm}
The following result is obtained in \cite{dlL08, SdlL12a} for $d=1$. It is shown that whenever the original functional
equation \eqref{FE'} has a variational structure, the additional parameter $\lambda$ vanishes when solving
\eqref{FE'}.
\begin{lemma}\label{vanishing lm 1}
 If  $U=\partial_{\alpha} V$, the solution
 $(h_{\ast}, \lambda_{\ast})$ of \eqref{FE'}
satisfies $\lambda_{\ast}=0$.
\end{lemma}
\noindent\textbf{Proof.}~Let $l_{\ast}= 1+\partial_{\alpha} h_{\ast}$ be a scalar function. We show that 
\begin{equation}\label{vanishing 1}
\langle l_{\ast}\cdot  (\Delta h_{\ast}+ \epsilon U(id+\alpha h_{\ast}))\rangle=0.
\end{equation}
Indeed, we have
\begin{equation*}
l_{\ast}(\theta)\cdot U(\theta+\alpha h_{\ast}(\theta))=\partial_{\alpha} V(\theta+\alpha h_{\ast}(\theta))\cdot (1+\partial_{\alpha} h_{\ast}(\theta))=\partial_{\alpha} [V(\theta+\alpha h_{\ast}(\theta)]
\end{equation*}
and thus $\langle l_{\ast}\cdot U(id+\alpha h_{\ast})\rangle=0$. In addition,
\begin{equation*}
\langle l_{\ast}\cdot \Delta h_{\ast}\rangle=\langle \Delta h_{\ast}\rangle+ \langle\partial_{\alpha} h_{\ast}\cdot \Delta h_{\ast}\rangle=\langle\partial_{\alpha} h_{\ast}\cdot \Delta h_{\ast}\rangle
\end{equation*}
since $\langle \Delta h_{\ast}\rangle=0$ by the definition of the difference operator $\Delta$ in \eqref{def Delta}.
Using the simple fact that 
\begin{equation*}
\int_{\TT^{m}} (T_{\omega_{j}\alpha} \partial_{\alpha} h_{\ast})(\theta)\cdot  h_{\ast}(\theta)~d\theta=
\int_{\TT^{m}}  \partial_{\alpha} h_{\ast}(\theta)\cdot (T_{-\omega_{j}\alpha}  h_{\ast})(\theta)~d\theta,
\end{equation*}
we conclude from
\begin{equation*}
\begin{aligned}
\int_{\TT^{m}} \partial_{\alpha} h_{\ast}(\theta)\cdot \Delta h_{\ast}~ d\theta=&-\int_{\TT^{m}} \partial_{\alpha} (\Delta h_{\ast})\cdot h_{\ast}(\theta)~d\theta
=-\int_{\TT^{m}} \Delta (\partial_{\alpha} h_{\ast})\cdot h_{\ast}(\theta)~d\theta\\
=&- \int_{\TT^{m}} \partial_{\alpha} h_{\ast}(\theta)\cdot \Delta h_{\ast}(\theta)~d\theta
\end{aligned}
\end{equation*}
that $\langle\partial_{\alpha} h_{\ast}\cdot \Delta h_{\ast}\rangle=0$ and hence \eqref{vanishing 1} holds.

Since $(h_{\ast},\lambda_{\ast})$ is a solution of \eqref{FE'}, it follows from $\langle l_{\ast}\rangle=1$
and \eqref{vanishing 1} that
\begin{equation*}
0=\langle l_{\ast}\scr{F}[h_{\ast}, \lambda_{\ast}]\rangle=\langle l_{\ast}\cdot  (\Delta h_{\ast}+ \epsilon U(id+\alpha h_{\ast}))\rangle+\langle l_{\ast} \lambda_{\ast}\rangle=\lambda_{\ast} \langle l_{\ast}\rangle=\lambda_{\ast},
\end{equation*}
which proves the lemma.\qed

\section{Proof of Theorem \ref{KAM}}\label{proof KAM}
We shall apply the Nash-Moser iteration method to solve nonlinear functional equation \eqref{FE'}.
Given an approximate solution $(h,\lambda)$ of \eqref{FE'}, we denote the error by 
\begin{equation}\label{def E}
E=\scr{F}[h, \lambda]\quad \textrm{with}~\|E\|_{\rho}< \bm{\varepsilon}.
\end{equation}
The Newton equation is 
$$
D\scr{F}[h, \lambda](\delta_{h}, \delta_{\lambda})=-\scr{F}[h, \lambda],
$$
or more precisely
\begin{equation}\label{Newton}
\Delta \delta_{h}+\epsilon \partial_{\alpha} U(id+ \alpha h)\cdot \delta_{h}+\delta_{\lambda}=- E.
\end{equation}

\subsection{Formulation of the quasi-Newton equation}
We shall employ the ideas in \cite{dlL08, SdlL12a} to modify the original Newton equation \eqref{Newton}
such that it can be simplified considerably. The key observation is that the variable coefficient
$\partial_{\alpha} U(\theta+\alpha h(\theta))$ in \eqref{Newton} can be obtained by
taking directional derivative of \eqref{def E}. Indeed, we see from \eqref{def F} and \eqref{def E}  that
\begin{equation}\label{E'}
\partial_{\alpha}E= \Delta (\partial_{\alpha} h)+ \epsilon \partial_{\alpha} U(\theta+\alpha h) \cdot (1+\partial_{\alpha} h).
\end{equation}
Let 
\begin{equation*}
l= 1+\partial_{\alpha} h,
\end{equation*}
and define the modified correction 
\begin{equation}\label{trans corr}
\widetilde{\delta_{h}}= l^{-1} \delta_{h}
\end{equation}
whenever $l^{-1}$ exists.
Notably, the directional derivatives of the hull functions (such as  $l$ and $\partial_{\alpha} U$)  are  scalars and therefore commute.

Solving $ \partial_{\alpha} U(\theta+\alpha h)$ from \eqref{E'} and substituting it into
\eqref{Newton}, we obtain the equation for $\widetilde{\delta_{h}}$ as follows
\begin{equation*}
\Delta (l\widetilde{\delta_{h}})- \Delta (\partial_{\alpha} h)\cdot \widetilde{\delta_{h}}+   \partial_{\alpha} E\cdot \widetilde{\delta_{h}}= -E-\delta_{\lambda}.
\end{equation*}
Note that, at each step of the Nash-Moser iterations, one would usually obtain $\|\widetilde{\delta_{h}}\|=O(\|E\|^{1-})$,
which implies that $\|\partial_{\alpha} E\cdot \widetilde{\delta_{h}}\|=O(\|E\|^{2-})$.
This suggests us to drop the quadratic error term $\partial_{\alpha} E\cdot \widetilde{\delta_{h}}$
for the moment and consider
\begin{equation}\label{Quasi Newton}
\Delta (l\widetilde{\delta_{h}})- \Delta l\cdot \widetilde{\delta_{h}}= -E-\delta_{\lambda},
\end{equation}
where we use the fact that $\Delta l= \Delta (\partial_{\alpha} h)$ since $\Delta 1=0$.

To transform the second order difference equation  \eqref{Quasi Newton} into 
first order, we define
\begin{equation*}
\nabla_{j} =T_{\omega_{j}\alpha} -id,\quad \nabla_{-j} =T_{-\omega_{j}\alpha} -id.
\end{equation*}
A direct computation yields that
\begin{equation*}
\begin{aligned}
\Delta (l\widetilde{\delta_{h}})-\Delta l\cdot \widetilde{\delta_{h}}=&\frac{1}{2d}\sum_{j=1}^{d} T_{\omega_{j}\alpha} l\cdot \nabla_{j} \widetilde{\delta_{h}}- T_{-\omega_{j}\alpha} (l\cdot \nabla_{j} \widetilde{\delta_{h}})\\
=&-l^{-1} \frac{1}{2d}\sum_{j=1}^{d}  T_{-\omega_{j}\alpha} (l\cdot T_{\omega_{j}\alpha} l\cdot \nabla_{j} \widetilde{\delta_{h}})
-l\cdot T_{\omega_{j}\alpha} l\cdot \nabla_{j} \widetilde{\delta_{h}}\\
=& -l^{-1} \frac{1}{2d}\sum_{j=1}^{d} \nabla_{-j} (l\cdot T_{\omega_{j}\alpha} l\cdot \nabla_{j} \widetilde{\delta_{h}}).
\end{aligned}
\end{equation*}
Then equation \eqref{Quasi Newton} becomes the following quasi-Newton equation
\begin{equation}\label{Quasi Newton 1}
L_{h} \widetilde{\delta_{h}}= l E+l \delta_{\lambda},
\end{equation}
where the linear operator $L_{h}$ is given by
\begin{equation*}
L_{h} v= \frac{1}{2d}\sum_{j=1}^{d} \nabla_{-j} (l\cdot T_{\omega_{j}\alpha} l\cdot \nabla_{j} v).
\end{equation*}

\begin{remark}
When $d=1$, the abstract operator $L_{h}$ is exactly the same to that in \cite{SdlL12a}, which is written
explicitly as two  first order difference equations that are uncoupled.
\end{remark}

\subsection{Analysis of the quasi-Newton equation}
In this part, we shall analyze the quasi-Newton equation \eqref{Quasi Newton 1}, especially
on the estimate of the operator $L_{h}^{-1}$. The main idea follows from \cite{dlL08}.
To begin with, we write
\begin{equation*}
L_{h}= S_{w} ~(S_{w}^{-1}~L_{h} S_{w}^{-1} )~S_{w}.
\end{equation*}
Let 
\begin{equation}\label{def b}
b_{j}(\theta)=l(\theta)\cdot l(\theta+\omega_{j}\alpha)-1=
 \partial_{\alpha} h(\theta)+\partial_{\alpha} h(\theta+\omega_{j}\alpha)+\partial_{\alpha} h(\theta)\cdot \partial_{\alpha} h(\theta+\omega_{j}\alpha).
\end{equation}
We have
\begin{equation*}
\begin{aligned}
S_{w}^{-1} L_{h} (S_{w}^{-1}v)= & S_{w}^{-1} \frac{1}{2d}\sum_{j=1}^{d} \nabla_{-j} \left((1+b_{j})\cdot \nabla_{j} (S_{w}^{-1}v)\right)\\
=&S_{w}^{-1} \frac{1}{2d}\sum_{j=1}^{d} \nabla_{-j}\nabla_{j}(S_{w}^{-1}v)+ S_{w}^{-1}\frac{1}{2d}\sum_{j=1}^{d} \nabla_{-j} (b_{j}\cdot \nabla_{j} (S_{w}^{-1}v))\\
=&id+S_{w}^{-1} \frac{1}{2d}\sum_{j=1}^{d} \nabla_{-j} (b_{j}\cdot \nabla_{j} (S_{w}^{-1}v)),
\end{aligned}
\end{equation*}
since 
\begin{equation*}
 \frac{1}{2d}\sum_{j=1}^{d} \nabla_{-j}\nabla_{j}v=\sum_{k\in\ZZ^{m}}  \frac{1}{2d}\sum_{j=1}^{d} 
 (e^{ 2\pi i \omega_{j}\langle k, \alpha\rangle}-1)(e^{- 2\pi i  \omega_{j}\langle k, \alpha\rangle}-1)
 \cdot \widehat{v}(k) e^{2\pi i \langle k,\theta\rangle}
\end{equation*}
and the definition of $S_{w}$ implies
\begin{equation*}
S_{w} ^{-1}  \frac{1}{2d}\sum_{j=1}^{d} \nabla_{-j}\nabla_{j} (S_{w}^{-1}v)=v.
\end{equation*}
This is the main motivation for our  definition of  $S_{w}$.

Let 
\begin{equation}\label{def K}
K_{h}=S_{w}^{-1} \left(\frac{1}{2d}\sum_{j=1}^{d} \nabla_{-j} (b_{j}\cdot \nabla_{j})\right) S_{w}^{-1}.
\end{equation}
For $v\in\scr{A}_{\rho}^{0}$, the average $\langle K_{h}v\rangle=0$ due
to Lemma \ref{Sw -1} and the fact that $\langle\nabla_{-j} v\rangle=0$.
The Hilbert matrix representation of $K_{h}$, denoted by $\scr{K}_{h}$, is given by
\begin{equation*}
\scr{K}_{h}(k,k')=\frac{1}{2d}\sum_{j=1}^{d}\frac{e^{-2\pi i \omega_{j}\langle k,\alpha\rangle}-1}{w(k)}\cdot 
\frac{e^{2\pi i \omega_{j}\langle k',\alpha\rangle}-1}{w(k')} \widehat{b}_{j}(k-k') .
\end{equation*}
It follows from \eqref{wei} that
\begin{equation*}
|\scr{K}_{h}(k,k')|\leq \sum_{1\leq j\leq d}|\widehat{b}_{j}(k-k')|,
\end{equation*}
and thus 
\begin{equation*}
\begin{aligned}
\|K_{h} v\|_{\rho}=&\sum_{k\in\ZZ^m}\left |\sum_{k'\in\ZZ^m} \scr{K}(k,k') \widehat{v}(k')\right| e^{2\pi \rho |k|}\\
\leq & \sum_{1\leq j\leq d} \sum_{k'\in\ZZ^{m}} \left(\sum_{k\in\ZZ^{m}}
|\widehat{b}_{j}(k-k')| e^{2\pi \rho |k-k'|}
\right) |\widehat{v}(k')| e^{2\pi \rho |k'|}\\
\leq & \left(\sum_{1\leq j\leq d} \|b_{j}\|_{\rho}\right)\cdot \|v\|_{\rho}.
\end{aligned}
\end{equation*}
In addition, if  $\sum_{1\leq j\leq d} \|b_{j}\|_{\rho}<\frac{1}{2}$,
we see from  Neumann series 
that $(id+K_{h})^{-1}\in\scr{L}(\scr{A}_{\rho}^{0},\scr{A}_{\rho}^{0})$ satisfying
$\|(id+K_{h})^{-1}\|_{\scr{L}(\scr{A}_{\rho}^{0},\scr{A}_{\rho}^{0})}\leq 2$.

In summary, we conclude the above results in the following proposition.

\begin{proposition}\label{Lh -1}
Let  $0<\sigma<\rho$ and $K_{h}$ be defined in \eqref{def K}. Then the linear operator $L_{h}$ can be written
as
\begin{equation*}
L_{h}= S_{w} (id+ K_{h}) S_{w},
\end{equation*}
where $K_{h}: \scr{A}_{\rho}^{0}\rightarrow \scr{A}_{\rho}^{0}$ satisfies 
\begin{equation*}
\|K_{h} \|_{\scr{L}(\scr{A}_{\rho}^{0}, \scr{A}_{\rho}^{0})}\leq \sum_{1\leq j\leq d} \|b_{j}\|_{\rho}.
\end{equation*} 

Moreover, if 
\begin{equation}\label{est b}
\sum_{1\leq j\leq d} \|b_{j}\|_{\rho}<\frac{1}{2},
\end{equation}
we have
\begin{equation*}
\|L_{h}^{-1}\|_{\scr{L}(\scr{A}_{\rho}^{0}, \scr{A}_{\rho-\sigma}^{0})}\leq C_{d,m,\tau}~\gamma^{-2}~ \sigma^{-2(\tau+m)} .
\end{equation*}
\end{proposition}

Now we consider the quasi-Newton equation \eqref{Quasi Newton 1}.
\begin{lemma}\label{sol quasi}
Let $0<\sigma<\rho$.
If 
\begin{equation}\label{cond h}
\max\left\{\|h\|_{\rho}, \|\partial_{\alpha}h\|_{\rho}\right\}< \frac{1}{10d},
\end{equation}
then  there is  a unique solution $(\widetilde{\delta_{h}^{0}}, \delta_{\lambda})$ of \eqref{Quasi Newton 1} 
 satisfying $\langle\widetilde{\delta_{h}^{0}}\rangle=0$ such that
\begin{equation*}
\begin{aligned}
&|\delta_{\lambda}|< 2\bm{\varepsilon},\\
&\|\widetilde{\delta_{h}^{0}}\|_{\rho-\sigma}\leq C_{d,m, \tau}\gamma^{-2} \sigma^{-2(\tau+m)} \bm{\varepsilon}.
\end{aligned}
\end{equation*}
\end{lemma}

\noindent\textbf{Proof.}~
If $\|\partial_{\alpha}h\|_{\rho}< \frac{1}{10d}$, then we see from \eqref{def b} that
the condition \eqref{est b} holds.
To ensure that $\langle l E+ l\delta_{\lambda}\rangle=0$, we choose
\begin{equation*}
\delta_{\lambda}= -\langle l E\rangle\cdot \langle l\rangle^{-1}=-\langle l E\rangle,
\end{equation*}
which satisfies
\begin{equation*}
|\delta_{\lambda}|\leq 2 \|E\|_{\rho}.
\end{equation*}
Then we have
\begin{equation*}
\widetilde{\delta_{h}^{0}}= L_{h}^{-1}(l E+ l\delta_{\lambda})
\end{equation*}
with $\langle\widetilde{\delta_{h}^{0}}\rangle=0 $
and Proposition \ref{Lh -1} implies
\begin{equation*}
\|\widetilde{\delta_{h}^{0}}\|_{\rho-\sigma}<C_{d, m,\tau}\gamma^{-2} \sigma^{-2(\tau+m)} \|E\|_{\rho}.
\end{equation*}
\qed

\begin{remark}
Since $\langle L_{h} v\rangle=0$,
the solution of  quasi-Newton equation \eqref{Quasi Newton 1} is not unique. 
This fact would be used to construct the correction $\delta_{h}$ of the Newton equation
such that $\langle\delta_{h}\rangle=0$.
\end{remark}

\subsection{
Newton iterative lemma}
In this part, we state the iterative lemma as follows.
Let
\begin{equation*}
\rho_{+}=\rho-2\sigma, \quad \bm{\varepsilon}_{+}= C\gamma^{-4}\sigma^{-4(\tau+m)-1} \bm{\varepsilon}^{2}.
\end{equation*}

\begin{lemma}
\label{Iter Lm}
Let $(h, \lambda)$ be an approximate solution of \eqref{FE'} with 
$\langle h\rangle=0$ such that  \eqref{cond h} and the estimate
$$\|E\|_{\rho}\leq \bm{\varepsilon}$$ hold.
Then there is  the correction $(\delta_{h}, \delta_{\lambda})$ satisfying 
$\langle \delta_{h}\rangle=0$ and 
\begin{equation*}
\max\Big\{\|\delta_{h}\|_{\rho_{+}}, \|\partial_{\alpha}\delta_{h}\|_{\rho_{+}}\Big\}\leq C_{d, m,\tau}\gamma^{-2}
\sigma^{-2(\tau+m)-1} \bm{\varepsilon},\quad 
|\delta_{\lambda}|\leq  2\bm{\varepsilon},
\end{equation*}
such that the new error $E_{+}=\scr{F}[h_{+}, \lambda_{+}]$
with 
\begin{equation*}
h_{+}= h+\delta_{h},\quad \lambda_{+}=\lambda+ \delta_{\lambda}
\end{equation*}
satisfies
\begin{equation*}
\|E_{+}\|_{\rho_{+}}\leq \bm{\varepsilon}_{+}.
\end{equation*}
\end{lemma}

\noindent\textbf{Proof.~}
Recalling the solution $(\widetilde{\delta_{h}^{0}},\delta_{\lambda})$ of
\eqref{Quasi Newton 1} in Lemma \ref{sol quasi}, we take
\begin{equation*}
\widetilde{\delta_{h}}= \widetilde{\delta_{h}^{0}}-\langle  l \widetilde{\delta_{h}^{0}}\rangle\quad \textrm{and}\quad \delta_{h}= l\widetilde{\delta_{h}}.
\end{equation*}
It follows that
$
\langle \delta_{h}\rangle= 0
$
since  $\langle l\rangle=1$.
For the Newton equation \eqref{Newton}, we take $(h_{+}, \lambda_{+})$  as the new approximate solution. 
Consider 
\begin{equation*}
\begin{aligned}
E_{+}=\scr{F}[h_{+}, \lambda_{+}]=\scr{F}[h, \lambda]+ D\scr{F}[h, \lambda](\delta_{h}, \delta_{\lambda})+\scr{R}[\delta_{h}],
\end{aligned}
\end{equation*}
where the Taylor remainder  is given by
\begin{equation*}
\begin{aligned}
\scr{R}[\delta_{h}]=&\epsilon \Big(
U(id+\alpha h+\alpha \delta_{h})- U(id+\alpha h)- \partial_{\alpha} U(id+\alpha h) \delta_{h}
\Big)\\
=&  \epsilon \int_{0}^{1}\int_{0}^{1} \partial_{\alpha}^{2} U(id+\alpha h+\mu\nu\alpha \delta_{h})\cdot \mu 
\delta_{h}^{2} ~d\mu d\nu.
\end{aligned}
\end{equation*}

The unsolved term in the Newton equation \eqref{Newton} is
$\partial_{\alpha} E\cdot \widetilde{\delta_{h}}$, which satisfies
\begin{equation*}
\|\partial_{\alpha} E\cdot \widetilde{\delta_{h}}\|_{\rho-\sigma}
\leq C_{d,m, \tau} \gamma^{-2} \sigma^{-2(\tau+m)-1} \|E\|_{\rho}^{2}.
\end{equation*}
The Taylor remainder satisfies
\begin{equation*}
\|\scr{R}[\delta_{h}]\|_{\rho-\sigma}\leq C_{d,m, \tau} \gamma^{-4} \sigma^{-4(\tau+m)} \|E\|_{\rho}^{2}.
\end{equation*}
As a result, we have
\begin{equation*}
\|E_{+}\|_{\rho-\sigma}\leq C_{d,m, \tau} \gamma^{-4} \sigma^{-4(\tau+m)-1} \bm{\varepsilon}^{2}.
\end{equation*}
\qed

\subsection{Convergence of the iteration and proof of Theorem \ref{KAM}}
The proof of the convergence of iterations is standard, which is provided
below for completeness.
For the initial step of the iteration, we take
\begin{equation*}
h_{0}=0,\quad \lambda_{0}=0.
\end{equation*}
As a result, we have
\begin{equation*}
\bm{\varepsilon}_{0}=\|E_{0}\|_{\rho_{0}}= \|\scr{F}[h_{0},\lambda_{0}]\|_{\rho_{0}}=\epsilon \| U\|_{\rho_{0}}.
\end{equation*}

At $n$-th step of the Nash-Moser iteration, we define the iterative constants as follows
\begin{equation*}
\sigma_{n} = \min\{1, 2^{-(n+1)} \sigma_{\ast}\},\quad \rho_{n+1}=\rho_{n}-2\sigma_{n},\quad \sigma_{\ast}=\rho_{0}/10.
\end{equation*}
and 
\begin{equation*}
\bm{\varepsilon}_{n+1}= C\gamma^{-4}\sigma_{n}^{-4(\tau+m)-1} \bm{\varepsilon}_{n}^{2}.
\end{equation*}
One readily sees that  $\rho_{\infty}= \lim_{n\rightarrow\infty} \rho_{n}\geq \rho_{0}/10 $.
Moreover, taking $$\tilde{\bm{\varepsilon}}_{n}=C \gamma^{-4}\sigma_{\ast}^{-(4(\tau+m)+1)} 2^{(4(\tau+m)+1) (n+2)}\bm{\varepsilon}_{n},$$ we have $\tilde{\bm{\varepsilon}}_{n+1}=\tilde{\bm{\varepsilon}}_{n}^{2}$
and thus 
$\tilde{\bm{\varepsilon}}_{n}= \tilde{\bm{\varepsilon}}_{0}^{2^{n}}$, which implies 
$\bm{\varepsilon}_{n}\rightarrow 0$ as long as 
\begin{equation}\label{small cond}
\tilde{\bm{\varepsilon}}_{0}= C\epsilon \gamma^{-4}\sigma_{\ast}^{-(4(\tau+m)+1)} 4^{4(\tau+m)+1}\|U\|_{\rho_{0}}< 1.
\end{equation}

Assume by induction that for $0\leq j\leq n$, there is an an approximate solution $(h_{j}, \lambda_{j})$   of \eqref{FE'} with 
$\langle h_{j}\rangle=0$ and $h_{j}= h_{j-1}+\delta_{h,j-1}, \lambda_{j}= \lambda_{j-1}+\delta_{\lambda,j-1}$ such that $$\|E_{j}\|_{\rho_{j}}=\|\scr{F}[h_{j}, \lambda_{j}]\|_{\rho_{j}}\leq \bm{\varepsilon}_{j},$$
and 
\begin{equation*}
\|\delta_{h, j-1}\|_{\rho_{j}}\leq C_{d,m, \tau}\gamma^{-2}
\sigma_{j-1}^{-2(\tau+m)-1} \bm{\varepsilon}_{j-1},\quad |\delta_{\lambda,j-1}|\leq 2\bm{\varepsilon}_{j-1},
\end{equation*}
in which we take $\delta_{h,-1}=0$ and $\delta_{\lambda,-1}=0$  for the consistence of notations.

If \eqref{small cond} holds, we also have
\begin{equation*}
\|h_{n}\|_{\rho_{n}}=\|h_{0}+\sum_{j=1}^{n} \delta_{h,j-1}\|_{\rho_{n}}\leq \sum_{j=1}^{\infty} C_{d,m, \tau}\gamma^{-2} \sigma_{j-1}^{-2(\tau+m)-1} \bm{\varepsilon}_{j-1}
\leq 2^{-(4(\tau+m)+1)}\gamma^{2} \sum_{j=0}^{\infty} \tilde{\bm{\varepsilon}}_{j}< \frac{1}{10d}.
\end{equation*}
Then, by Lemma \ref{Iter Lm}, the induction statement hold for $j=n+1$.
Similar arguments show that $\lambda_{n}$ converges to some constant $\lambda_{\infty}$
and $\lim_{n\rightarrow\infty} \|h_{n}-h_{\infty}\|_{\rho_{\infty}}=0$ for some function $h_{\infty}\in \scr{A}^{0}_{\rho_{\infty}}$.
Moreover, the limit function $(h_{\infty}, \lambda_{\infty})$ is a solution of \eqref{FE'}.

Finally, combining with the vanishing Lemma \ref{vanishing lm 1}, we complete the proof of Theorem \ref{KAM}.

\section{Discussion}\label{sect dis}

In this section, we briefly explain the possibility of applying the Craig-Wayne-Bourgain (CWB) method
to construct the quasi-periodic equilibrium of \eqref{FE},
especially the connection of the problem with some existing results in the
spectral theory of lattice Schr\"{o}dinger operators.
We only outline some ideas of applying
the CWB method to the F-K model and refer to
\cite{Bou94_IMRN, Bou97_MRL, Bou98_Ann, Bou05_Green,CW93, HSSY24} for more
technical details on nonlinear problems.
We shall also compare the CWB method with the classical KAM method in
\cite{dlL08, Koz83,SdlL12a} as well as the present paper.

Following the formulation and notations in Section \ref{sect FE}, we consider \eqref{FE}
in the form
\begin{equation}\label{QP NL 1}
\Delta h(\theta)+\epsilon\partial_\alpha V(\theta+\alpha h(\theta))=0,
\qquad \theta\in\TT^m,
\end{equation}
where $\omega\in\RR^d$, $\alpha\in\RR^m$ and $\Delta$ is defined in \eqref{def Delta}.
Applying the Nash-Moser iteration, we need to solve the Newton equation with variable coefficient
\[
\Delta\delta_h+\epsilon\Phi\delta_h=-E,
\qquad \Phi(\theta)=\partial_\alpha^2 V(\theta+\alpha h(\theta)),
\]
where $E$ is the left-hand side of \eqref{QP NL 1} evaluated at the approximate solution $h$.
Passing into Fourier coefficients, we obtain
\begin{equation}\label{abs T}
\mathcal T\widehat{\delta_h}=(\mathcal D+\epsilon\mathcal S)\widehat{\delta_h}=-\widehat E,
\end{equation}
where $\mathcal D$ is diagonal and $\mathcal S$ is a Toeplitz matrix:
\[
\mathcal D(k)=\frac1d\sum_{j=1}^d
\bigl[\cos 2\pi(\omega_j\langle k,\alpha\rangle)-1\bigr],
\qquad \mathcal S(k,k')=\widehat\Phi(k-k'),\quad k,k'\in\ZZ^m.
\]
If $\Phi\in\scr A_\rho$, then
$|\mathcal S(k,k')|\leq\|\Phi\|_\rho e^{-2\pi\rho|k-k'|}$,
where $|\cdot|$ is the $\ell^1$ norm used above.
Thus $\mathcal T$ is a long-range quasi-periodic lattice operator on $\ZZ^m$.
For real approximate solutions its finite-volume restrictions are self-adjoint,
and their inverse norms are the reciprocals of the distances from zero to their spectra,
whenever the inverses exist. The Newton equation therefore concerns Green's functions at energy zero.

At each step of the Newton iteration, we would apply the truncation technique to solve $\delta_h$,
with the support of $\widehat{\delta_h}$ contained in $[-N,N]^m$, from
\[
\mathcal T_N\widehat{\delta_h}=-\widehat E_N,
\]
where $\mathcal T_N$ and $\widehat E_N$ are the restrictions to $[-N,N]^m\cap\ZZ^m$.
The main task is to establish the Green's function estimates for the original operator:
\begin{equation}\label{Green}
\begin{aligned}
&\|\mathcal T_N^{-1}\|_{\ell^2\to\ell^2}<A^{(\log N)^{C_1}},\\
&|\mathcal T_N^{-1}(k,k')|<\exp(-c|k-k'|),
\quad \textrm{for }|k-k'|>(\log N)^{C_2},
\end{aligned}
\end{equation}
for $C_2>C_1>0$, $c>0$ and a large constant $A>1$, on the frequencies retained in the iteration.

To obtain \eqref{Green}, we introduce an auxiliary parameter $\sigma\in\TT^d$ and define
\[
\mathcal T^\sigma=\mathcal D^\sigma+\epsilon\mathcal S,
\qquad
\mathcal D^\sigma(k)=\frac1d\sum_{j=1}^d
\bigl[\cos 2\pi(\sigma_j+\omega_j\langle k,\alpha\rangle)-1\bigr].
\]
Here $\mathcal T^0=\mathcal T$, and the family satisfies the covariance relation
\[
\mathcal T^\sigma(k+p,k'+p)
=\mathcal T^{\sigma+\langle p,\alpha\rangle\omega}(k,k'),
\qquad p\in\ZZ^m.
\]
As long as we establish a large deviation theorem for this auxiliary family,
with a small exceptional set of phases, the next step is to eliminate the auxiliary
parameter through a further parameter-exclusion argument and obtain \eqref{Green}
for the original operator outside a small set of frequencies.
This step uses the covariance relation and  the decomposition theorem for semialgebraic sets,
see \cite{HSSY24} for example.

In proving the large deviation theorem, a key ingredient needed is a bound on the number
of bad local blocks inside a larger  block.
Let $\mathcal B_M\subset\TT^d$ contain the phases at which the scale-$M$ Green's function estimates fail.
By covariance, the centers to be counted are
\[
\bigl\{p\in[-N,N]^m\cap\ZZ^m:
\sigma+\langle p,\alpha\rangle\omega\in\mathcal B_M\bigr\},\qquad M<N.
\]
A suitable bound on this number, together with initial-scale estimates, matrix-valued Cartan
estimates and multiscale coupling, allows one to establish the large deviation estimate
at the next scale, see \cite{Bou05_Green,HSSY24} for example.
The required orbit-counting estimates take particularly simple forms in the following two cases.

When $d=1$ and $m\geq2$ (i.e., the quasi-periodic equilibrium problem for one dimensional F-K model), the auxiliary phase belongs to $\TT$.
Writing $\beta=\omega\alpha$, condition \eqref{SDC} gives
$\operatorname{dist}(\langle q,\beta\rangle,\ZZ)\geq\gamma|q|^{-\tau}$ for $q\ne0$.
Hence the phases corresponding to distinct $p\in[-N,N]^m\cap\ZZ^m$ are separated
by at least $\gamma(2mN)^{-\tau}$.
If $\mathcal B_M$ is covered by $J$ intervals of total length at most $\eta$, each interval
of length $a$ contains at most $1+a\gamma^{-1}(2mN)^\tau$ such phases. Consequently,
\[
\#\bigl\{p\in[-N,N]^m\cap\ZZ^m:
\sigma+\langle p,\beta\rangle\in\mathcal B_M\bigr\}
\leq J+\eta\gamma^{-1}(2mN)^\tau.
\]
For a semialgebraic cover, $J$ is controlled by its degree.
Controlling the number and lengths of the intervals, and choosing compatible scales,
therefore gives the required bound of the number of  bad blocks.

When $m=1$ and $d\geq2$ (i.e., the periodic equilibrium problem for the generalized F-K model), we normalize $\alpha=1$ and the phase orbit is
$\sigma+k\omega\in\TT^d$.
Suppose that a semialgebraic cover $S_M\supset\mathcal B_M$ has degree $B$ and measure less than $\eta$,
and that $\omega$ satisfies a Diophantine condition
$\operatorname{dist}(\langle q,\omega\rangle,\ZZ)\geq\gamma_0|q|^{-\nu}$
for all $q\in\ZZ^d\setminus\{0\}$, with fixed $\gamma_0,\nu>0$.
Then \cite[Corollary~9.7]{Bou05_Green} gives, uniformly in $\sigma$,
\[
\#\bigl\{1\leq k\leq N:\sigma+k\omega\in S_M\bigr\}<N^{1-\delta}
\]
for some $\delta=\delta(\omega)>0$, provided $N$ is sufficiently large and
$\log B\ll\log N<\log(1/\eta)$.
The same estimate, up to a constant, applies to $-N\leq k\leq N$.
Here the Diophantine assumption is an additional requirement: \eqref{SDC} only controls
$\max_j\operatorname{dist}(k\omega_j,\ZZ)$ in this case and does not imply this bound.
These two counting mechanisms explain the connection with the multiscale analysis of
quasi-periodic lattice Schr\"{o}dinger operators: arithmetic controls how often local
Green's function estimates fail along the phase orbit.

As we can see from the discussion above, an advantage of the CWB
approach is that it does not rely on the particular factorization
\eqref{Lh trans}.
It studies the linearized operator through its finite-volume
Green's functions and allows for resonant local blocks, provided
that their number and interaction can be controlled.
This offers a possible approach to extensions of the equilibrium
problem in which the factorization is no longer available, while
retaining the connection with the spectral theory of quasi-periodic
lattice operators.

On the other hand, compared with condition \eqref{SDC}, the frequency
conditions required by the proposed CWB approach can be more
restrictive at two different stages.
First, the bad-block counting bound may require stronger arithmetic
assumptions, as in the application of
\cite[Corollary~9.7]{Bou05_Green} when $m=1$ and $d\geq2$.
Second, passing from the large deviation theorem for the auxiliary
family to the Green's function estimates at $\sigma=0$ involves
a further frequency exclusion through the decomposition theorem
for semialgebraic sets.
In comparison, the present method uses the special structure of
the equilibrium equation.
The factorization \eqref{Lh trans} separates the small-divisor
estimates for $S_w^{-1}$, which follow from \eqref{SDC}, from
the inversion of $id+K_h$, which follows from the Neumann series
under the smallness assumptions of Theorem \ref{KAM}.
The construction therefore applies throughout the prescribed
simultaneous Diophantine class, without either of these additional
frequency restrictions.
\section*{Acknowledgment}
The authors would like to thank Prof. Rafael de la Llave for bringing our attention
to the paper \cite{Koz83}. X.H. would thank Prof. Xifeng Su for several fruitful discussions.
 X.H. supported by NSFC
(Grant NO.12371167 and NO. 12471156). Z.Y. supported by NSFC (Grant NO. 12371189).

\bibliographystyle{plain}
\def\cprime{$'$} \def\cprime{$'$} \def\cprime{$'$}
\bibliography{reference}

\end{document}